\documentclass[12pt]{amsart}

\usepackage{amssymb,amsbsy,amsmath,amsfonts,amssymb,amscd}

\input xy
\xyoption{all}

\newcommand\sC{{\mathcal C}}

\newcommand\sG{{\mathcal G}}

\newcommand\sN{{\mathcal N}}
\newcommand\sK{{\mathcal K}}
\newcommand\sS{{\mathcal S}}

\newcommand\la{\lambda}

\newcommand\Ga{\Gamma}
\newcommand\Lam{\Lambda}
\newcommand\De{\Delta}

\newcommand{\CC}{\ensuremath{\mathbb{C}}}

\newcommand{\PP}{\ensuremath{\mathbb{P}}}

\newcommand{\ra}{\ensuremath{\rightarrow}}

\def\eea{\end{eqnarray*}}
\def\bea{\begin{eqnarray*}}

\newcommand\dual{\mathrel{\raise3pt\hbox{$\underline{\mathrm{\thinspace d
\thinspace}}$}}}

\newcommand\QED{\ifhmode\unskip\nobreak\fi\quad {\rm Q.E.D.}} 
\newcommand\qe{\ifhmode\unskip\nobreak\fi\quad $\Box$}       

\def\BOX{\hfill\lower.5\baselineskip\hbox{$\Box$}}

\newtheorem{theo}{Theorem}[section]
\newtheorem{remarkk}[theo]{Remark}
\newenvironment{rem}{\begin{remarkk}\rm}{\end{remarkk}}

\newtheorem{defin}[theo]{Definition}

\newtheorem{prop}[theo] {Proposition}

\newtheorem{lemma}[theo]{Lemma}
\newtheorem{example}[theo]{Example}

\newtheorem{problem}[theo]{Problem}
\newtheorem{rema}{Remark}[section]

\def\La{\langle\langle}
\def\Ra{\rangle\rangle}

\newcommand{\Proof}{{\it Proof. }}
\begin{document}

\title[Caustics  ]{Caustics of plane  curves, their birationality and matrix projections.}
\author{Fabrizio Catanese}


\thanks{The present work took place in the realm of the DFG
Forschergruppe 790 "Classification of algebraic surfaces and compact
complex manifolds". }

\date{\today}

\maketitle

{\em  Dedicated to Klaus Hulek
on the occasion of his $60$-th birthday.}

\date{\today}

\maketitle

\section{Introduction and  setup}

Given a plane curve $C$ and a point $S$, a source of light (which could also lie at infinity, as the sun),
the light rays $L_P$ originating in $S$, and hitting the curve $C$ in a point $P$, are reflected by the curve, and
the caustic $\sC$ of $C$ is the envelope of the family of reflected rays $\Lam_P$.

Our first  Theorem  \ref{caustic} says that the correspondence between  the curve $C$ and  the caustic
curve $\sC$ is birational, i.e., it is generically one to one,  if the light source point $S$ is chosen to be a general point.

We learnt about this problem in  \cite{JP1}, to which we refer for an account of the history of the theory of caustics and for
references to the earlier works of von Tschirnhausen, Quetelet, Dandelin, Chasles, and more modern ones (as \cite{bgg1}, \cite{bgg2}).

Our methods are   from   algebraic geometry,   so we got interested in a generalization of this result, in which the special form
of a certain curve $D$ plays no role: we achieve this goal in  Theorem \ref{MP}.

 Let us   now describe  the mathematical set up  for the description of caustics.

Let $\PP^2 = \PP^2_{\CC}$ and let $ C \subset \PP^2$ be a plane irreducible algebraic curve, whose normalization
shall be denoted by $C'$.

 Choose an orthogonality structure in the plane, i.e. two points,
called classically  the cyclic points, and let $\PP^1_{\infty}$ be the line (`at infinity') joining them.
The two cyclic points determine a unique  involution $\iota$ on  $\PP^1_{\infty}$ for which the cyclic points are fixed,
hence an involution, called orthogonality, on  the pencils of lines passing through a given point of the affine plane $\PP^2 \setminus \PP^1_{\infty}$.

Without loss of generality, we choose  appropriate projective coordinates such that 

$$ \iota : (x,y,0)  \mapsto (-y,x,0), \\ Fix (\iota) = \{ (1, \pm \sqrt{-1}, 0)\}.$$ 

Let $S \in \PP^2$ be a light source point, and to each point $P \in \PP^2 \setminus \{S\}$ associate the line $L_P : = \overline {PS}$.
In the case where $P \in C$, we define $\Lam_P$, the reflected light ray, as the element of the pencil of lines
through $P$ determined by the condition that the cross ratio
$$ CR ( N_P, T_P, L_P, \Lam_P) = -1,$$
ensuring the existence of a symmetry with centre $P$ leaving the tangent line $T_P$ to $C$ at $P$ 
and the normal line $N_P : = \iota (T_P)$ fixed, and exchanging the incoming  light ray $L_P$ with the reflected
light ray $\Lam_P$.

We thus obtain a rational map of the algebraic curve $C$  to the dual projective plane:

$$ \Lam : C \dashrightarrow (\PP^2)^{\vee} .$$
\begin{defin}  The Caustic $\sC$ of $C$ is defined as the envelope of the family of lines $\{\Lam_P\}$: in other words,  setting  $\Ga : = \Lam (C)$, $\sC  = \Ga^{\vee}$.
\end{defin}

\begin{rem} since the biduality map $\Ga \dashrightarrow  \Ga^{\vee}$ is birational (cf. \cite{Walker}, pages 151-152), the map $C \dashrightarrow \sC$ is birational iff
$\Lam : C \dashrightarrow \Ga$ is birational. Moreover, by the biduality theorem,  the class of the caustic $\sC$  is the degree of $\Ga$, and the degree of $\sC$
is the class of $\Ga$. 
\end{rem} 

We shall quickly see  in the next section the basic  calculations  which give  the class of  $\sC$, i.e. the  degree of $\Ga$, 
 in the case where $C$
and $S$ are general (more precise Pl\"ucker type formulae which show how the singularities of the
curve $C$ and the special position of $S$ make these numbers decrease are to be found in \cite{JP1} and \cite{JP2}).

In section 3 we show our fist  result,  that $\Lam$ is birational onto its image  for general choice of
the source point $S$, if $C$ is not a line (in this case $\Ga$ is a line, and  the caustic is  a point). The next  section recalls a  well known lemma about lines contained in
the determinantal variety $\De$ which is the secant variety of the Veronese surface $V$.

This lemma plays a crucial role in the proof of our main result,  which
 says  the following (see Theorem \ref{MP} for more details):

\begin{theo}

Let $     D \subset \PP :  = \PP (Sym^2 (\CC^3))$ be a curve.

Then, for general $ S \in \PP^2$,
the  projection $\pi_S :  \PP = \PP (Sym^2 (\CC^3)) \dashrightarrow \PP^2$ given by
$ \pi_S (B) : = B S$ has the property that its restriction to $D$, $\pi_S | _D$, is  birational 
onto its image, unless (and this is indeed an exception) $D$ is a curve contained in a  plane $ \De (S') = \{ B | \ B S' =  0 \}$ (contained
in the determinantal hypersurface  $ \De = \{ B | \det (B) = 0\}$) and $D$  is not a line.

\end{theo}

This result  suggests the investigation of a more general situation concerning the birationality
of linear projections given by matrix multiplications.

\begin{problem}
Given a linear space $\PP$ of matrices $ B$,
and a linear space $\PP'$  of matrices $ S$, consider the matrix multiplication
$ \pi_S (B) = BS$.  For which algebraic  subvarieties $D \subset \PP$ is the restriction
of the projection $ \pi_S | _D$ birational onto its image for a general choice of $S \in \PP'$ ?
\end{problem}

\section{Equations in coordinates}

Let $f (x_0, x_1,x_2)=0$ be the equation of $C$ in the appropriate system of homogeneous coordinates,  let $d : = deg (f)$,
and let $F:= ((f_0(x), f_1(x))$ be the first part of the gradient of $f$.
For a point $x = (x_0, x_1,x_2)$ we define $$ ( F,x ): =  f_0(x) x_0 +  f_1(x)x_1, \ \ 
\{F\wedge x\}:  =  f_0(x) x_1 -  f_1(x)x_0.$$

Then the  tangent line $T_P$ at a point $ P $ with coordinates $ x$ is the transpose of  the row vector $(f_0(x), f_1(x),f_2(x))$.
 
The normal  line $N_P$ is orthogonal to the tangent line, hence it has the form 
$N_P =  \ ^t(- f_1(x), f_0(x),f_3(x))$, and the condition that $ P \in N_P$ forces the unknown rational function $f_3(x)$
to fulfill $  - f_1(x) x_0 +  f_0(x) x_1 + f_3(x) x_2 \equiv 0$, thus

$^tN_P$ is  the row vector 
$$ ^tN_P = (- x_ 2 f_1(x) , x_ 2 f_0(x)  , - \{F\wedge x\}).$$

We find now the line $L_P$ as the line in the pencil spanned by $T_P$ and $N_P$ passing through $S$:
as such the line $L_P$ is a column  vector which is a linear combination  $\la T_P + \mu N_P$;
the condition that $ S \in L_P$ then determines $\la = - ^tN_P \cdot S , \mu = ^tT_P \cdot S$,
 where $S$ is the transpose
of the vector $(s_0, s_1,s_2)$.

Hence  we get 

$$ L_P (S) = A(P) S ,   \ \ \ A (P) : =  - T_P \ ^tN_P + N_P \ ^tT_P,$$
in particular the matrix $A(P)$ is skew symmetric.

To obtain the reflected ray $\Lam(P)$ it is sufficient, by definition, to change the sign of $\la$, and we get therefore:

$$ \Lam_P (S) = B(P) S ,   \ \ \  \ B  (P) : =   T_P \ ^tN_P + N_P \  ^tT_P. $$

\begin{rema}\label{B}

1) The matrices $A(P)$ and $B(P)$ are functions 
which are defined for all general points $P$ of the plane.

2) The matrix $B(P)$ is  symmetric and has rank at most two, since its image is generated by $N_P$ and $T_P$;
moreover we have 
$$ B(P) P = 0 , A (P) P = 0, \ \  \forall P \in C.  $$ 

3) Assume that $C$ is not a line passing through a cyclic point: then  the matrix $B(P)$ has precisely rank two on the 
non empty open set  where  $f_1^2 + f_0^2 \neq 0$ and $x_2 \neq 0$; 
the former  condition clearly holds
for a general point $ P\in C$, otherwise the dual curve of $C$ would be contained in a line  $ y_0 = \pm \sqrt{-1} y_1$.

4) The entries of the matrix $B (x) $ are given by polynomials of degree $2d-1$.

\end{rema}

By the preceding  remark follows easily the classical theorem asserting that
\begin{theo}
The class of the caustic, i.e., the degree  of $\Ga$, equals  $ d (2d-1)$, for a general curve $C$ and a general choice of $S$. 

\end{theo}

In fact $C$ has degree $d$, and  $ B(x) S$ is given by polynomials
of degree $2d-1$ in $x$, which have no base points on a general curve $C$.

\section{Birationality of the caustic map}

\begin{theo}\label{caustic}
 If $C$ is not a line, then the caustic map $ C \dashrightarrow \sC$ is birational, for general choice of $S$. 

\end{theo}

\Proof

As already remarked, the caustic map is birational iff the map $\Lam : C \dashrightarrow  \Ga$ is birational.
Observe that $\Lam$ defines  a morhism $ C' \ra \Ga$ which we also denote by $\Lam$.

The matrix $B$, whose  entries are polynomials of  degree $2d-1$, yields
a map 
$$   B : C' \ra  D \subset \PP^5 = \PP (Sym^2 (\CC^3)).$$

\begin{lemma}
$  B : C'  \ra  D: = \Phi (C)$ is birational.

\end{lemma}

\Proof
It suffices to recall  remark \ref{B}: for a general point $P \in C$, 
 $B(P)$ has rank exactly two, and $ B (P) P = 0$. Hence  $ P =  \ker (B(P))$, and  the  matrix $B(P)$
determines the point $ P \in \PP^2$.

\qed

We have now a projection $ \PP (Sym^2 (\CC^3)) \dashrightarrow \PP^2$ given by

$$ \pi_S (B) : = B S.$$

Consider  the linear subspace $$ W :  = \{ B |  B_{0,0} + B_{1,1} = 0 \}.$$

We observe preliminarily that the curve $D$ is contained in the linear subspace 
$W $
since, setting for convenience $f_i : = f_i (x)$,  the matrix $B(x)$ has  the following entries:
$$ B_{0,0} =  - 2 x_2 f_0 f_1 , \ \ B_{1,1} =  2 x_2 f_0 f_1 .$$ 

Then our main result follows from the next assertion, that,  for a general choice of $S \in \PP^2$,
 the projection  $\pi_S$ yields a birational map
of $D$ onto $\Ga := \pi_S (D)$.

In order to prove this, we set up  the following notation:

 $$\De_S  :  = \{ B |  B S = 0 \}, \ \ \    \De : = \{ B | \det (B) = 0 \} = \cup_S \De_S .  $$ 
 
 Observe that $\De$ is the secant variety of the Veronese surface
 $$ V : = \{ B  | \ rank (B) = 1\}.$$

Observe that the curve $D$ is  contained in the linear subspace $ W  ,$
is contained in $\De$ but not contained in the Veronese surface $V$.

We are working inside the subspace $W$, and we observe first of all that the centre of the projection $\pi_S$
restricted to $W$ 
is the linear space 
$$ W_S : = \De_S \cap W.$$

Observe moreover  that   $ \De \cap W =  \cup_S W_S.$ 

Now, the projection $\pi_S$ is  not birational on $D$ if and only if, for a general  $B \in D$,
there exists another $B' \in D, \ B \neq B'$, such that the chord (i.e., secant line) $ B * B'$ intersects 
$W_S$ in a point $B''$ (observe that the general point $B \in D$ does not lie in the line $W_S$). 

There are two possible cases:

{\bf Case I}:

$B''$ is independent of the point $B \in D$.

{\bf Case II}:

$B''$ moves as a rational function  of the point $B \in D$, hence the points $B''$ sweep the line $W_S$.

\begin{lemma}
The assumption that case I holds for each $S \in \PP^2$ leads to a contradiction.
\end{lemma}
{\em Proof of the Lemma.}
Under our assumption, for each $S $ there is a point $B'' (S) $ such that infinitely many chords of $D$ meet $W_S$ in $B'' (S) $.

Let us see what happens if we specialize  $S$ to be a general point $P \in C$. 

The first alternative is 

I-1) $B'' (P) = B (P)$: in this case, for each point $B^1  \in D$ there is $B^2  \in D$ such that $ B(P), B^1, B^2$
are collinear. Since this happens for each choice of $B(P), B^1$, every secant is a trisecant,
hence , by the well known trisecant lemma (cf. \cite{ACGH}, page 110), $D$ is a plane curve of order at least three.

Take now a general $S \in \PP^2$: since $B''(S)$ is on a secant to $D$, $B''(S)$ belongs to the secant variety $\Sigma$ 
of $D$ (here a plane
$\Pi$), but we claim that it is not in $D$.
In fact, if there were a point $P \in C'$ such that  $B''(S) = B(P)$, then  $B(P) S = 0$ contradicting that  $S$ is a general
point.
Hence we obtain that the plane $\Pi$  intersects $\De$ in a bigger
locus than $D$: since $\De$ is a cubic hypersurface, it follows that $\Pi \subset \De$. 

By proposition \ref{LS} it follows that either there is a point $S'$ such that $ S' \in \ker (B), \forall B \in \Pi$,
or there is a line $ L \in \PP^2$ such that   $  \ker (B) \in L, \forall B \in \Pi$: both cases imply that the curve $C$ must be contained in a line,
a contradiction.

The second alternative is 

I-2 )  $B'' : =  B'' (P) \neq B (P)$.
Then there is a point $B' \in D$ (possibly infinitely near) such that $B'$ is a  linear combination of $B''$ and $B : = B(P)$.
 
  However, since 
$ B P = 0, B'' P = 0, $ and  $ B \neq B''$, then also for their linear combination $B'$ we have $B' P = 0$.
The consequence is, since $ B' P = B' P' = 0$, that $B'$
has rank one.  Therefore, if $B'$ is not infinitely near,  $B'$ cannot be a general point of $D$, hence $B'$ is independent of
$P$: but then $C \subset \ker (B') $, and since we assume that $C$ is not a line, we obtain $B' = 0$, a contradiction.

If $P'$ is  infinitely near to the  point $P \in C$, i.e.,  $ P, P'$ span the tangent line to $C$ at $P$,
and $ B, B'$ span the tangent line to $D$ at $B = B(P)$, we work over the ring of tangent vectors $\CC[\epsilon] / (\epsilon^2)$, and we observe that 
$$  ( B + \epsilon B') ( P +  \epsilon P') = 0 \Rightarrow  B P' = 0.$$

For $P \in C$ general this is a contradiction, since $ B P' = 0, BP = 0$ imply that $B = B(P)$ has rank one.

\qed

\begin{lemma}
The assumption that case II  holds for general  $S \in \PP^2$ leads to a contradiction.
\end{lemma}
{\em Proof of the Lemma.}
As we already observed, for general $S$, $B''$ moves as a rational function  of the point $B \in D$, hence the points $B''$ sweep the line 
$W_S$. Therefore the line $W_S$ is contained in the secant variety $\Sigma$ of the curve $D$.
As this happens for general $S$, and  $ \De \cap W =  \cup_S W_S,$
it follows that the threefold $ \De \cap W $ is contained in the secant variety $\Sigma$.

Since $\Sigma$ is irreducible, and has dimension at most three, it follows that we have equality
$$ \De \cap W  = \Sigma .$$

We conclude that, for $P_1, P_2$ general points of $C$, the line joining  $B(P_1)$ and $ B( P_2)$ is contained in $\De$.

By proposition \ref{LS}, and since  $ker (B(P_1) )= P_1 , ker  (B(P_2) )= P_2$,  we have that the  matrices in  the pencil 
$ \la_1 B(P_1) + \la_2 B(P_2) $ send the span of $P_1, P_2$ to its orthogonal subspace.

This condition is equivalent to 

$$ ^t P_1 (B(P_2) )  P_1 = 0 \ \forall P_1, P_2  \in C$$ 
($ ^t P_2 (B(P_1) )  P_2 = 0$ follows in fact since $P_1, P_2$ are general).

Fix now a  general point $P_2$: then we have a quadratic equation for $C$, hence $C$ is contained in a conic.

A little bit more of attention: the matrix $B(P_2)$ has rank two, hence the quadratic equation defines a reducible
conic, and, $C$ being irreducible, $C$ is a line, a contradiction.

\qed

\section{Linear subspaces contained in the determinantal cubic $\De : = \{ B | \det (B) =0 \}$}

\begin{prop}\label{LS}
Let $\la B_0 + \mu B_1$ be a line contained in the determinantal hypersurface $\De$ of the projective space of symmetric
$ 3 \times 3$ matrices.

Then the line is contained in a maximal projective subspace contained in $\De$, which is either of the type

$$ \De_S : =  \{ B |  B S = 0\},$$
for some $S \in \PP^2$, or of the type  

$$ \De (L)   : =  \{ B |  B L \subset L^{\perp}  \} = \{ B |  B _{|  L} \equiv 0   \} ,$$
for some line  $L \subset \PP^2$.
\end{prop}

\Proof

A pencil of reducible conics either has at most one (non infinitely near) base point $ S \in \PP^2$, or it has a line $L$ as fixed component.

In the first case the pencil is $ \subset \De_S$, in the second case it is contained in 
the subspace $ \De (L)$ consisting of the conics of the form $ L + L'$,
where $L'$ is an arbitrary line in the plane.

\qed 

\begin{rem}
Even if the result above  follows right away  from the classification of pencils of conics, it is useful to  recall the arguments which will be used in the sequel.

For instance, we observe that the hyperplane sections of the Veronese surface $V$ are smooth conics, hence no line is contained in $V$.
\end{rem}

\section{Birationality of certain matrix projections of curves}

In this final section we want to show the validity of a much more general statement:

\begin{theo}\label{MP}

Let $     D \subset \PP :  = \PP (Sym^2 (\CC^3))$ be a curve and 
$   B : C' \ra  D \subset \PP $ be its normalization. 

Then, for general $ S \in \PP^2$,
the  projection $\pi_S :  \PP = \PP (Sym^2 (\CC^3)) \dashrightarrow \PP^2$ given by
$ \pi_S (B) : = B S$ has the property that its restriction to $D$, $\pi_S | _D$ is  birational 
onto its image, unless $D$ is a curve contained in a plane $ \De (S') $ and is not a line.

In the latter case, each  projection $\pi_S | _D$ has as image the line $ (S')^{\perp}$
and is not birational.
\end{theo}

\Proof 

Let $ \sG : = Gr (1, \PP)$ be the Grassmann variety of lines $\Lam \subset \PP$: $\sG$ has dimension 8.

Define, for $S \in \PP^2$, $\sG_S: = \{ \Lam \in \sG | \Lam \cap \De_S \neq \emptyset \}.$ Indeed,
these 6-dimensional submanifolds of $\sG$ are the fibres of the second projection of the
incidence correspondence
$$I \subset \sG \times \PP^2, \ I : = \{ (\Lam, S) |  \Lam \cap \De_S \neq \emptyset \} .$$
In   turn $I$ is the projection of the correspondence 
$$J \subset \sG \times \De \times \PP^2, \ J : = \{ (\Lam, B, S) | B \in  \Lam , \   B S = 0 \} .$$

Recall further that $\De \setminus V$ has a fibre bundle structure 
$$ \sK : \De \setminus V \ra \PP^2 $$
such that $ \sK (B) : = ker (B)$, and with fibre over $S$ equal to $\De_S \setminus V$.

\begin{rem}\label{ImB}
(1) Observe that for matrices $ B \in V$ we can write them in the form $ B = x \  ^tx$,
for a suitable vector $x$, and in this case $ \ker (B) = x^{\perp}$, $Im(B) = \La x \Ra$.

(2) In any case, since the matrices $B$ are symmetric, we have always 
$$ Im (B) = \ker (B)^{\perp}. $$
\end{rem}

Consider now the fibres of $ I \ra \sG$: for a general line $\Lam$, its fibre $\sS (\Lam )$ is

\begin{enumerate}
\item
if $\Lam \cap \De \neq \Lam , \Lam \cap \De \subset \De \setminus V$, then $\sS (\Lam )$ consists of at most three points;
\item
if $\Lam \cap \De \neq \Lam , \ ( \Lam \cap  V) \neq \emptyset$, then $\sS (\Lam )$ consists of a line  $x^{\perp}$
and at most one further point;
\item
if $ \Lam \subset \De$ is of the form   $ \Lam \subset \De_S$, then $\sS (\Lam )$ consists of one or two lines containing $S$;
\item
if $ \Lam \subset \De$ is of the form $ \Lam \subset  \De (L)$, $\sS (\Lam )$ consists of  the line $L$.

Since, if $ \Lam \subset  \De (L)$, the conics in $ \Lam $ consist of $L$ plus a line $L'$ moving in the pencil of lines
through a given point $P$.

\end{enumerate}

We let $$U \subset  \sG \times \PP , U : = \{ (\Lam, B) | B \in \Lam \}$$
be the universal tautological $\PP^1$-bundle, and we denote by $ p : U \ra \PP$
the second projection.

Recall now that the secant variety $\Sigma$ of $D$ is defined as follows:
we have a rational map $\psi :  C' \times C' \dashrightarrow \sG$ associating to the pair
$(s,t)$ the line $B(s) * B(t)$ joining the two image points $B(s), B(t)$.

Then one denotes by $U'$ the pull back of the universal bundle, and defines
$\Sigma$ as the closure of the image $p(U')$.

The condition that for each $S \in \PP^2$ the projection $\pi_S$ is not birational on $D$ means that,
if $Y$ is the closure of the image of $\psi$, then $ Y \cap \sG_S$ has positive dimension.

This implies that the correspondence 
$$ I_D : = \{ (\Lam_y, S) | y \in Y, \  \Lam_y \cap \De_S \neq \emptyset \} \subset Y \times \PP^2$$
has dimension at least three and surjects onto $\PP^2$.

Projecting $I_D$ on the irreducible surface $Y$, we obtain that all the fibres have positive dimension,
and we infer that each secant line $\Lam_y$ has a fibre $\sS (\Lam_y)$ of positive dimension.

There are two alternatives: 

(i) a general secant  $\Lam_y$ is not contained in $ \De$, but intersects the Veronese surface $V$.

(ii) each secant line $\Lam_y \subset \De$.

Step I) : the theorem holds true if $D \subset V$.

Proof of step I.

In this case any element of $D$ is of the form $B(t) = x(t) ^t x(t)$, and 
$$ \pi_S (B(t)) = x (t) [ ^t x(t) S]  = (x(t), S)  x(t)  = x(t) .$$ Hence, for each $S$, the projection
$\pi_S$ is the inverse of the isomorphism
$$\phi:  x \in \PP^2 \ra V,  \phi (x) = x \  ^t x.$$
\qed

We may therefore assume in the sequel that $D$ is not contained in $V$.

Step II) :  the theorem holds in case (i).

Proof of step II.

Observe preliminarly that, in case (i), $ D \not\subset \De$; else we could take two smooth points $B_1, B_1 \in  D \cap (\De \setminus V)$,
and the secant line $B_1 * B_2$ could not fulfill (i).

Choose then a point $ B_0 \in D, B_0  \in \PP \setminus \De$, hence w.l.o.g. we may assume that $B_0$ is the identity
matrix $I$.

Since any other point $B(t) \in D$ is on the line joining $B_0$ with a point $x(t) ^t x(t) \in V$,
we may write locally around a point of $C'$
$$  B(t) = I + \xi (t) ^t \xi (t) ,$$
where $\xi (t)$ is a vector valued holomorphic function.

Now, for each $s, t ,$ the secant line $ B(t) * B(s)$ meets the Veronese surface $V$.

Since $B(t)$ cannot have rank equal to 1, there exists $\la$ such that

$$ \la B(t) + B(s) = \la ( I + \xi (t) ^t \xi (t) ) + (I + \xi (s) ^t \xi (s)) $$
has rank equal to 1, i.e., 
$$ K_{\la}: = \ker ( \la B(t) + B(s) ) = \{ v |  [\la ( I + \xi (t) ^t \xi (t) ) + (I + \xi (s) ^t \xi (s))] v = 0 \} = $$
$$ \{ v |  (\la + 1) v + \la  \xi (t) ( \xi (t) , v ) +  \xi (s) (\xi (s), v )= 0 \} $$
has dimension 2.

Let us now make the assumption:

(**) two general points $ \xi (t), \xi(s)$ are linearly independent.

The above formula shows however that, under  assumption (**),
it must be that $v$ is a linear combination of  $ \xi (t), \xi(s)$. This is clear if $\la + 1 \neq 0$,
otherwise $v$ is orthogonal to the span of  $ \xi (t), \xi(s)$, contradicting that the kernel has dimension 2.

Hence $ K_{\la} = \La   \xi (t), \xi(s) \Ra$ and the condition that $\xi (t) \in K_{\la} $ yields
$$ (\la + 1) \xi (t)  + \la  \xi (t) ( \xi (t) , \xi (t)  ) +  \xi (s) (\xi (s), \xi (t)  )= 0$$
and implies 
$$ (***) \  \forall s,t \  (\xi (s), \xi (t)  )= 0.$$

(***) says that  $ K_{\la} = \La   \xi (t), \xi(s) \Ra$ is an isotropic subspace, which can have at most dimension 1.

Hence assumption (**) is contradicted, and we conclude that it must be:
$$ (****)  \  \xi (t) = f(t) u, $$
where $u$ is an isotropic vector and $ f(t)$  is a scalar function.

Even if this situation can indeed occur, we are done since in this case the matrix in $V$ is unique, 
$ u  ^t u$, each secant $\Lam_y$ contains $ u  ^t u$, hence $\sS (\Lam_y) = u^{\perp} \cup T_y$
where $T_y$ is a finite set.  Therefore, for general $S$, the fibre  $ \{ y | \Lam_y \cap \De_S \neq \emptyset \}$
is a finite set.

\qed

Step III: the theorem holds true  in case (ii).

Proof of Step III.

Consider the general secant line $ \Lam_y$. We have two treat two distinct cases.

Case (3) :  $ \Lam_y \subset \De$ is of the form  $ \Lam_y \subset \De_S$ (then $\sS (\Lam_y )$ consists of one or two lines containing $S$).

Case (4):  $ \Lam_y \subset \De$ is of the form $\Lam_y \subset \De (L)$ ( then $\sS (\Lam_y )$ consists of  the line $L$).

\medskip 
In case (3), this means that two general matrices $B(s), B(t)$ have a common kernel $S(s,t)$. Since the general
matrix $B(t)$ is in $\De \setminus V$, its rank equals 2 and $ S(s,t) = S(t) \  \forall s$.

Hence the curve $D$ is contained in a plane $\De_S$. In this case however $ Im B(t) \subset S^{\perp}$
and every projection $\pi_{S'} (B(t) = B (t) S' $ lands in the line $S^{\perp}$, so that the projection cannot be birational, unless our curve $D$ is a line.

\medskip

In case (4) for two general matrices $B(s), B(t)$  there exists a line $L = L (s,t)$ such that $B(s), B(t ) \in \De (L)$. 

Since two such general matrices have rank equal to 2, and $ B(t) L  \subset L^{\perp},  B(s) L  \subset L^{\perp},$
if $ v(t) \in  \ker B(t)$  it follows that $v(t) \in L$ ( since $ \ker B(t) \cap L \neq \emptyset$). Therefore, if 
$  B(t) \neq B(s)$, then $ L (t,s) = \La  v(t), v(s)\Ra $.

However, the above conditions  $ B(t) L  \subset L^{\perp},  B(s) L  \subset L^{\perp}$
are then equivalent to 
$$   ( B(t)  v(s), v(s) ) = ^t v(s) B(t) v(s) = 0, \ \forall t,s .$$

Fixing $t$ this is a quadratic equation in $v(s)$, but, since the curve $D$ is irreducible, and $B(t)$ has rank equal to 2,
we see that the vectors $ v(s)$ belong to a line. Therefore the line $L = L (s,t)$ is independent of $s,t$
and the conclusion is that the curve $D$ is contained in the plane $\De (L)$.

In suitable coordinates for $\PP^2$, we may assume that $L = \La  e_2, e_3 \Ra $ and $L^{\perp} = \La  e_1 \Ra$.

Choosing then $S = e_1$, we obtain an isomorphic projection, since for a matrix

\[ B = 
  \left(
  \begin{matrix}
    a &    b  &   c \\
     b & 0 &  0  \\
     c & 0  & 0
  \end{matrix} \right)
\] 
we have 
\[  B (e_1 ) =
  \left(
  \begin{matrix}
    a \\
     b  \\
     c 
  \end{matrix}\right).
\]

\qed

\qed

\begin{rem}
The referee suggested some arguments to simplify the proofs.

For  Theorem \ref{caustic}, this is the proposal:

a) Firstly, in the case of the caustic, the curve $D$ parametrizes the reducible conics of the form $ T_P + N_P$,
where $T_P$ is the tangent to the curve $C$ at $P$, and $N_P$ is the normal.

If $S$ is a general point in $\PP^2$, then  the degree of $D$ equals the number of such conics passing through $S$,
hence, if $\nu$ is the degree of the curve $\sN$ of normal lines, $\mu$ is the degree of $C' \ra \sN$, then 
$$  deg (D) = deg (C^{\vee}) + \mu \nu.$$

The above formula shows that $deg (D) \geq 4$ if $C$ is not a line.

 In fact, then $deg (C^{\vee}) \geq 2$, while in general  $\nu \geq 1$ (the normal $N_P$ contains $P$).  But,  
if $\nu = 1$, then the dual curve of $\sN$, the evolute, is a point, so $C$ is a circle, but in this case $\mu = 2$.

b) Therefore, if one shows that $D$ is contained in a plane $\pi$, then the plane $\pi$ is contained in the cubic hypersurface
$\De$, hence we can apply Proposition \ref{LS}.

c) In turn, to show that $D$ is a plane curve, it is necessary and sufficient to show that  two general tangent lines
to $D$ meet, which  follows  if one  proves that:

d)  for each secant line there is a cone over $D$ and with
vertex a point $B''$, such that the secant line passes through $B''$ 

(since then the two tangent lines are coplanar).

In case I), d) follows since then, for each general $S$, there is a point $B''(S)$ such that a curve of secants passes through $B''(S)$,
and we get a cone over $D$ with vertex $B''(S)$. Varying $S$, the point $B''(S)$ must vary, since  $B''(S) S = 0$; hence the cone varies,
and we get that for each secant d) holds true.

In case II), as we have shown, the secant variety of $D$ equals  $ W \cap \De$, which is the secant variety of the the rational normal quartic
$ W \cap V$: but the singular locus of the secant variety of $ W \cap V$ equals $ W \cap V$ and contains $D$, hence $ W \cap V = D$,
a contradiction.

The argument suggested for Theorem \ref{MP} requires some delicate verification, so we do not sketch it here.

\end{rem}

{\bf Acknowledgement:}

I  would like to thank  Alfrederic Josse and Francoise P\`ene for stimulating email correspondence 
and for spotting a mistake in my
first naive attempt to prove birationality of the caustic map for general source $S$, thus  pushing me to  find the  proof 
 of Theorem \ref{caustic}, which I  announced to them   in an e-mail on December 20, 2012.

At the moment of writing up the references for the present article, I became aware, by searching on 
the arXiv, that they have written an independent and different proof of birationality of the caustic map
for general source,  in \cite{JP3}.

Thanks to the referee for helpful comments, and for sketching  alternative arguments, which are reproduced in the remark above.

\medskip
\noindent {\bf Author' s  Address:}\\
\noindent Fabrizio Catanese, \\ Lehrstuhl Mathematik VIII, Mathematisches
Institut der \\Universit\"at Bayreuth\\ NW II,  Universit\"atsstr. 30,
95447
Bayreuth, Germany.\\
            e-mail:          fabrizio.catanese@uni-bayreuth.de\\

\end{document}